\documentstyle[12pt]{article}
\begin{document}
\newcommand{\QED}{\hspace{0.2in}\vrule width 6pt height 6pt depth 0pt
\vspace{0.1in}}
\newcommand{\Cal}{\cal}
\newcommand{\Proof}{{\em Proof} \hspace{0.2in}}

\newcommand{\tag}[1]{\eqno(#1)}
\newenvironment{remark}{{\bf Remark} \hspace{0.2in}}{\\}
\newtheorem{theorem}{Theorem}[section]
\newtheorem{lemma}[theorem]{Lemma}
\newtheorem{proposition}[theorem]{Proposition}
\newtheorem{corollary}[theorem]{Corollary}
\newtheorem{problem}[theorem]{Problem}
\newtheorem{definition}[theorem]{Definition}

\begin{center}
{\LARGE Quantization of quadratic Poisson brackets\\
\medskip
on a polynomial algebra of three variables}\\
\bigskip
{\large J. Donin}\\
Department of Mathematics and Computer Science\\
Bar-Ilan University,
52900 Ramat-Gan, Israel\\
\smallskip
{\large L. Makar-Limanov}\\
Department of Mathematics and Computer Science\\
Bar-Ilan University,
52900 Ramat-Gan, Israel\\
and Department of Mathematics\\
Wayne State University,
Detroit 48202, USA
\end{center}
\date{}

\def\W{\wedge}
\def\t{\bigotimes}

\newcommand{\ff}{\varphi}
\def\d{\partial}
\def\k{\mbox{\bf k}}
\def\h{\hbar}
\def\ad{\mbox{ad}\,}
\def\wh{\widehat{W}(E)}
\def\wt{W(E)}
\def\ah{\widehat{A}}
\def\at{\widetilde{A}}
\def\fb{\bar{f}}

\def\ns{\mbox{nonsingular }}
\def\na{\nabla}
\def\df{\mbox{Der}_f(A)}

\def\Ah{A_{\h}}
\def\hb{\h}
\def\Ch{\k[[\hb]]}
\begin{abstract}
Poisson brackets (P.b) are the natural initial terms for the
deformation quantization of commutative algebras.
There is an open problem whether any Poisson bracket on the polynomial
algebra of $n$ variables can be quantized. It is known
(Poincare-Birkhoff-Witt theorem) that any linear P.b. for all $n$
can be quantized. On the other hand, it is easy to show that
in case $n=2$ any P.b. is quantizable as well.

Quadratic P.b. appear as the unitial terms for the quantization
of polynomial algebras as quadratic algebras.
The problem of the quantization of quadratic P.b. is also open.
In the paper we show that in  case $n=3$ any quadratic P.b can be
quantized. Moreover, the quantization is given as the quotient algebra
of tensor algebra of three variables by relations which are similar
to those in the Poincare-Birkhoff-Witt theorem. The proof uses a
classification of all quadratic Poisson brackets of three variables,
which we also give in the paper. In Appendix we give explicit algebraic
constructions of the quantized algebras appeared here to show that they
are related to algebras of global dimension three considered by M.Artin,
W.Schelter, J.Tate, and M.Van Den Bergh from a different point of view.
\end{abstract}

\def\I{{\cal I}}
\def\J{{\cal J}}
\def\t{\otimes}
\section{Deformations of quadratic algebras}
Let $A$ be an associative algebra with unit over a field $\k$ of
characteristic zero.
We will consider deformations of $A$ over the algebra of formal
power series $\Ch$ in a variable $\hb$.

By a deformation of $A$ we mean an algebra $\Ah$ over $\Ch$
which is isomorphic to $A[[\hb]]=A\hat{\t}_{\k}\Ch$ as a $\Ch$-module and
$\Ah/\h\Ah=A$ (the symbol $\hat{\t}$ denotes the tensor product
completed in the $\hb$-adic topology). We will also denote $A$ as $A_{0}$.

If $A'_{\hb}$ is another deformation of $A$, we call the deformations
$A_{\hb}$ and $A'_{\hb}$ equivalent if there exists a $\Ch$-algebra
isomorphism $A_{\hb}\to A'_{\hb}$ which induces the identity
automorphism of $A_0$.

Let $T=T(V)$ be a tensor algebra over a finite-dimensional vector
space $V$ over a field $\k$ and let $A$ be a quotient algebra of $T$, i.e.,
$A=T/\I$ where $\I$ is an ideal in $T$.
It is easy to see that if $A_{\h}$ is a deformation of $A$, then
$A_{\h}=T[[\h]]/\I_{\h}$ where $\I_{\h}$ is an ideal in $T[[\h]]$
such that $\I=\I_{\h}/\h\I_{\h}$.

Consider $T$ as a graded algebra, $T=\oplus_kT^k$, and let $A$ also be
a graded algebra, i.e., $\I$ is a graded ideal. Denote by $\I^k$ the
$k$-th homogeneous component of $\I$, $\I^k=\I\cap T^k$.
If $A_\h$ is a homogeneous deformation of $A$, then the ideal $\I_{\h}$
is homogeneous, and  for flatness of $A_\h$ it is necessary that
all the components $\I^k_{\h}$ be splitting  submodules in
$V^{\otimes k}$ for all $k$. This means that each $\I^k_\h$ has
a complementary submodule $\J^k_\h$, i.e., $V^{\t k}[[\h]]=\I^k_\h\oplus
\J^k_\h$.

Now consider the case of quadratic algebras. Recall that an algebra
$A$ is called quadratic if $A=T/\I$ and $\I$ is a homogeneous ideal
generated by $\I^2$. From now on we use the notations $\I^2=I$ and
$\I^k=I^k$ and perceive each $I^k$ as a subspace in $V^{\otimes k}$.
So, a quadratic algebra is defined by a pair $(V,I)$ where $I$
is a subspace in $V\otimes V$.
Define the subspaces $I^k_i,\,i=1,...,k-1$, in the linear space
$V^{\otimes k}$ as $I^k_i=V\otimes\cdots\otimes I\otimes \cdots\otimes V$
where $I$ is in the $i$-th position. It is clear that
the component $I^k$ is equal to the sum $\sum_{i=1}^{k-1}I^k_i$
of the subspaces $I^k_i$ of the linear space $V^{\otimes k}$.

A deformation of the quadratic algebra $(V,I)$ is defined by
a splitting $\k[[\h]]$-submodule $I_{\h}$ of $V^{\t 2}[[\h]]$ such that
all the submodules $I^k_{\h}=\sum_{i=1}^{k-1}I^k_{\h i}$ are splitting.
We shall write $V$ instead of $V[[\h]]$ when it does not lead to
misunderstanding. We call the dimension of the vector space $I^k_0$ over
$\k$ the rank
of the submodule $I^k_{\h}$ at the zero point. We call the dimension
of the vector space $I^k_{\h}\t_{\k[[\h]]}\k\{\h\}$ over the quotient
field $\k\{\h\}$ of $\k[[\h]]$ the rank
of the submodule $I^k_{\h}$ at the general point.
Note that $I^k_\h$ is splitting if and only if its ranks at the
zero point
and at the general point coincides.

In some cases a deformation
of a quadratic algebra is flat if $I^3_\h$ is a splitting submodule.
We call a quadratic algebra $(V,I)$  weakly Koszul if
it satisfies the following distributiveness condition:
$$I_1^k\cap(I_2^k+I_3^k+\cdots +I_{k-1}^k)=I_1^k\cap I_2^k+
I_1^k\cap (I_3^k+\cdots +I_{k-1}^k) \tag{1}$$
for all $k$.

This definition is motivated by the fact that a quadratic algebra
will be Koszul if and only if for all $k$ the lattice of subspaces
$I^k_i,i=1,...,k-1$, of
$V^{\otimes k}$ is distributive under the both operations:
taking  the sum and the intersections of subspaces (see \cite{Ba}).

The following result essentially follows from \cite{D} (see also \cite{PP}).
\begin{proposition}
\label{2.1a}
Let $I_\h$ be a deformation of a subspace $I\subset V^{\otimes 2}$.
In case $(V,I)$ is a weakly Koszul algebra, all $I^k_{\h}$ will be splitting
$\k[[\h]]$-submodules if $I^3_{\h}=V\otimes I_{\h}+I_{\h}\otimes V$ is
a splitting $\k[[\h]]$-submodule.
\end{proposition}
\Proof
The assertion is valid for $k=3$ by hypothesis of the proposition.
Suppose it is valid for all $m<k$. To prove that the $\k[[\h]]$-submodule
$I^k_\h=I_\h\t V^{\t(k-2)}+V\t I_\h^{k-1}$ is splitting, it is sufficient
to show that the module
$J_\h=I_\h\t V^{\t(k-2)}\cap V\t I_\h^{k-1}$ is splitting.
We will prove that the rank of $J_\h$ at the general point coincide
with the rank of this submodule at the zero point, which implies
that $J_\h$
is a splitting submodule.
We have the inclusion
$$I_\h\t V^{\t(k-2)}\cap
(V\t I_\h\t V^{\t(k-3)} + V^{\t 2}\t I_\h^{k-2})\supseteq
I_\h^{(3)}\t V^{\t(k-3)}+I_\h\t I_\h^{k-2}, \tag{2}$$
where $I^{(3)}_\h$ denotes $I_\h\t V\cap V\t I_\h$.

Now observe that the submodule
$(V\t I_\h\t V^{\t(k-3)} + V^{\t 2}\t I_\h^{k-2})=V\t I_\h^{k-1}$
is splitting by the induction hypothesis, and the submodules
$I_\h\t V^{\t(k-2)}$, $I_\h^{(3)}\t V^{\t(k-3)}$, and $I_\h\t I_\h^{k-2}$
are splitting as well. At the point $\h=0$ the modules in (2) are equal
because of the weakly Koszul property.
But it is clear that the rank
of the intersection of two splitting $\k[[\h]]$-submodules at the general
point has to be
not more than the rank of their intersection at the zero point,
while the rank of the sum of two splitting $\k[[\h]]$-submodules at the general
point has to be not less than the rank of their sum at the zero point.
Therefore, the ranks of both modules in (2)
 at the general point must be equal to
the rank of these modules at the zero point. Hence, the modules coincide
and $J_\h$ is a splitting submodule.
The proposition is proved.

{\bf Remark.}
Let $I_\h$ be an algebraic family of subspaces of $V\t V$, where point
$\h$ runs over a connected algebraic variety $M$,
such that the dimension of subspaces
$I_\h\t V+V\t I_\h$ of $V^{\t 3}$ does not change.
The arguments above show that if $I_o$ satisfies (1) at a
point $o\in M$, then for a fixed $k>3$ the dimension of $I_\h^k$
is constant for $\h\in M\setminus S_k$, where $S_k$ is an algebraic subset of
a smaller dimension. Therefore, for almost all points of $M$ (the complement to
$\cup S_k$) the dimension of all $I_\h^k$ does not change.\\

Given a quadratic algebra $(V,I)$, one can define the dual quadratic
algebra $(V^*, I^*)$ where $V^*$ is the dual space to $V$, and $I^*$
consists of the elements $\lambda\in V^*\otimes V^*$ such that
$\lambda(I)=0$.
It follows from duality that in order for $(V^*, I^*)$ to be weakly
Koszul, $I$ must satisfy the dual to (1) distributiveness condition:
$$I_1^k+(I_2^k\cap I_3^k\cap\cdots \cap I_{k-1}^k)=(I_1^k+ I_2^k)\cap
(I_1^k+ (I_3^k\cap\cdots \cap I_{k-1}^k) \tag{3}$$
for all $k$.
In particular, if the quadratic algebra $(V, I)$ is Koszul, then
its dual quadratic algebra $(V^*, I^*)$ will be Koszul as well.

If $(V,I_{\h})$ is a deformation of the algebra
$(V,I)$, one can build the dual $\k[[\h]]$-submodule $I^*_{\h}$
in $V^*[[\h]]$.

\begin{proposition}
\label{Pr2.1}
Let $(V,I_{\h})$ be a deformation of a Koszul quadratic algebra $(V,I)$.
Then the pair $(V^*,I^*_{\h})$ gives a deformation of the dual
quadratic algebra $(V^*,I^*)$.
\end{proposition}
\Proof The dual quadratic algebra $(V^*,I^*)$ also is a Koszul one,
so that it suffices to check that the $\k[[\h]]$-module
$I^{*3}_{\h}=V^*\otimes I^*_{\h}+I^*_{\h}\otimes V^*$ is splitting.
But $I^{*3}$ is dual to $V\otimes I_{\h}\cap I_{\h}\otimes V$,
which is a splitting $\k[[\h]]$-submodule because
$V\otimes I_{\h}+I_{\h}\otimes V$ is splitting.
So $I^{*3}$ is a splitting submodule as well, which completes the proof.\\

Let $x_1,...,x_n$ be a basis in the space $V$, and $k\leq n$.
We will perceive the first $k$ variables $x_1,...,x_k$ as even
and the other $x_{k+1},...,x_n$ as odd. Then the tensor algebra
$T(V)$ can be regarded as a superalgebra. The super-commutative
polynomial superalgebra $P_{k,l}$ ($l=n-k$) is the quotient
of $T(V)$ by the ideal generated by the elements
$$x_i\otimes x_j-(-1)^{x_i x_j}x_j\otimes x_i.$$
(Here $(-1)^z=1$ if $z$ is even and $(-1)^z=-1$ if $z$ is odd.)
These elements determine a subspace $I_{l,k}$ in $V\otimes V$.
Let $I_{k,l}$ be a subspace determined by the elements
$$x_i\otimes x_j+(-1)^{x_i x_j}x_j\otimes x_i.$$
The subspaces $I_{k,l}$ and $I_{l,k}$ are complemented in $V\otimes
V$, so that the natural projection $\pi:T(V)\to P_{k,l}$ gives
an isomorphism $\pi_2:I_{k,l}\to P^2_{k,l}$. Here $P^2_{k,l}$ denotes
the component of the degree two in $P_{k,l}$.
Further we set $I=I_{l,k}$ and $J=I_{k,l}$.
A deformation $I_{\h}$ of the subspace $I$ can be given by a
family of linear operators $B_{\h}:I[[\h]]\to J[[\h]],\ B_{\h}=\h B^1+\h^2
B^2+\cdots$, such that $I_{\h}$ is a graph of the operator $B_{\h}$.
Define a bilinear form $b$ on $P_{k,l}$ in the following way:
If $\bar{x}_i$, $\bar{x}_j$ are the images of $x_i$, $x_j$ under
the projection $\pi$, we put $b(\bar{x}_i, \bar{x}_j) =
\pi B^1(x_i\otimes x_j - (-1)^{x_ix_j}x_j\otimes x_i)$, and on others
elements from $P_{k,l}$ the form $b$ is extended by the Leibniz rule
$$b(uv,w)=ub(v,w)+(-1)^{uv}vb(u,w).$$
It is clear that in such a way we obtain a one-one correspondence
between the operators $B^1:I\to J$ and the skew-symmetric bilinear
forms $b:P_{k,l}\otimes P_{k,l}\to P_{k,l}$ satisfying the Leibniz
rule and such that $b(\bar{x}_i,\bar{x}_j)$ are homogeneous quadratic
forms,
$b(\bar{x}_i,\bar{x}_j)=\sum_{pq}b_{ij}^{pq}\bar{x}_p\bar{x}_q$.
We call such a quadratic bracket
Poisson if it satisfies the Jacobi identity
$$(-1)^{uw}b(b(u,v),w)+(-1)^{vu}b(b(v,w),u)+(-1)^{wv}b(b(w,u),v).$$
One can show in the standard way that if the family of operators
$B_{\h}$ determines a deformation of quadratic algebra
$P_{k,l}=(V,I)$, then
the form $b$ corresponding to $B^1$ is a quadratic Poisson bracket.
It is easy to see that the dual to $P_{k,l}$ quadratic algebra
can be identified with $P_{l,k}$, the supercommutative polynomial
superalgebra with $l$ even variables $x^1,...,x^l$ and $k$ odd
variables $x^{l+1},...,x^n$, $n=l+k$.
The quadratic Poisson bracket $b$ on $P_{k,l}$ defines a quadratic
Poisson bracket $\tilde{b}$ on $P_{l,k}$ in the following way:
$\tilde{b}(x^i,x^j)=\sum_{pq}b^{ij}_{pq}x^px^q$.
This bracket corresponds obviously to the dual deformation of the dual
algebra $P_{l,k}$.

A quadratic Poisson bracket $b$ on $P_{k,l}$ is called quantizable if
there exists a deformation of $P_{k,l}$ with $b$ as an initial term.
Any deformation of $P_{k,l}$ with initial term $b$ will be called a
quantization
(or realization) of $b$.
Since $P_{k,l}$ is a Koszul algebra, it follows from Proposition
\ref{Pr2.1} that
\begin{proposition}
\label{Pr2.2}
If the bracket $b$
is quantizable on $P_{k,l}$ then the dual bracket $\tilde{b}$ is
quantizable on $P_{l,k}$.
\end{proposition}

In the next sections we describe all quadratic Poisson brackets
on the polynomial algebra of three even variables and
show that any such bracket can be quantized, i.e. has a realization.
Using Proposition \ref{Pr2.2} we can conclude that any quadratic
Poisson bracket on the Grassmann algebra of three variables is
quantizable as well.

\section{Classification of quadratic Poisson brackets in
three-dimensional case}

\label{S3}
Let $b(f,g)$ be a Poisson bracket defined on
a polynomial ring
$R = \k[x_1,...,x_n]$
for
which all $y_{i,j} = b(x_i,x_j)$ are homogeneous quadratic forms.
Let $P(b)$ be the matrix with elements
$p_{i,j} =
\frac{\partial}{\partial x_j} (\sum_k \frac{\partial y_{i,k}}{\partial x_k})$.
Let $A$ be a matrix from $GL(n,\k)$. Matrix $A$ defines
an automorphism $a$ of $R$ by $x_i^a = \sum_j a_{i,j} x_j$. It is clear
that the same Poisson bracket may be given in terms of $\{x_i^a\}$.
Let us denote this Poisson bracket by $b^a$ and let $P^a = P(b^a)$.

\begin{lemma} $P^a = APA^{-1}$.
\end{lemma}
\Proof It is sufficient to check that the lemma is true for
elementary automorphisms
$x^a_i = x_i$ for $i \neq j$ and $x^a_j = x_j + cx_k$
($SL(n)$ is generated by these transformations)
and for diagonal automorphisms
$x^a_i = x_i$ for $i \neq j$ and $x^a_j = cx_j$.
This may be done by a straightforward computation.\\
(Another way to prove this lemma is to consider $b$ as an element of
$V^* \otimes V^* \otimes V \otimes V$. Then $P$ corresponds to the
contraction of this tensor.)

{\bf Remark.} Since $y_{i,j} = - y_{j,i}$ one can see that $trace P = 0$.\\

Let $b$ be a quadratic Poisson bracket and let $P$ be its corresponding
matrix. Let $r = rank P$.
One obtains a natural classification of quadratic Poisson brackets
according to the rank and Jordan form of $P$.\\
Let us denote by $v_i = \sum_k \frac{\partial y_{i,k}}{\partial x_k}$.
For $n = 3$ the Jacobi
identity is equivalent to $y_{1,2}v_3 + y_{2,3}v_1 + y_{3,1}v_2 = 0$.\\

{\bf a)} $r = 0$. In this case, all $v_i = 0$ which implies that there
exists a cubic form $f$ such that
$y_{1,2} = \frac{\partial f}{\partial x_3}$,
$y_{2,3} = \frac{\partial f}{\partial x_1}$,
and $y_{3,1} = \frac{\partial f}{\partial x_2}$.\\

{\bf b)} $r = 1$. In this case, there exists a linear automorphism after which
$v_1 = v_2 = 0$ and $v_3 = x_1$. Therefore, $y_{1,2} = 0$. But then
since $v_1 = v_2 = 0$,  derivatives
$\frac{\partial y_{1,3}}{\partial x_3} =
\frac{\partial y_{2,3}}{\partial x_3} = 0$
and $v_3 = x_1 = \frac{\partial y_{3,1}}{\partial x_1}
+ \frac{\partial y_{3,2}}{\partial x_2}$.
This means that there exists a cubic form $f$ such that
$y_{1,2} = \frac{\partial f}{\partial x_3} = 0$,
$y_{2,3} = \frac{\partial f}{\partial x_1} - x_1x_2$,
$y_{3,1} = \frac{\partial f}{\partial x_2}$.\\

{\bf c)} $r = 2$. Here we have two subcases:\\
{\bf ca)} $v_1 = 0$, $v_2 = -\lambda x_2$, and $v_3 = \lambda x_3$, or \\
{\bf cb)} $v_1 = 0$, $v_2 = x_1$, and $v_3 = x_2$.  \\
ca) There exists a cubic form $f$ such that
$y_{1,2} = \frac{\partial f}{\partial x_3}$,
$y_{2,3} = \frac{\partial f}{\partial x_1} - \lambda x_2x_3$,
$y_{3,1} = \frac{\partial f}{\partial x_2}$.
The Jacobi identity then gives
$x_3\frac{\partial f}{\partial x_3} - x_2\frac{\partial f}{\partial x_2} = 0$,
and it is easy to see that $f =  2c_1x_1x_2x_3 + c_2x_1^3$.\\
cb) There exists a cubic form $f$ such that
$y_{1,2} = \frac{\partial f}{\partial x_3}$,
$y_{2,3} = \frac{\partial f}{\partial x_1}
+ x_1x_3 - \frac{1}{2}x_2^2$,
$y_{3,1} = \frac{\partial f}{\partial x_2}$.
The Jacobi identity then gives
$x_2\frac{\partial f}{\partial x_3} +
x_1\frac{\partial f}{\partial x_2}   = 0$,
and here $f = -2c_1x_3x_1^2 + c_1x_1x_2^2 + c_2x_1^3$.\\

{\bf d)} $r = 3$. Here we have three subcases:\\
{\bf da)} $v_1 = \lambda _1x_1$, $v_2 = \lambda _2x_2$,
and $v_3 = \lambda _3x_3$,
where  $\lambda _1 + \lambda _2 + \lambda _3 = 0$ and
$\lambda _1$, $\lambda _2$, and $\lambda _3$ are pairwise different (and
different from zero).\\
{\bf db)} $v_1 = \lambda x_1$, $v_2 = \lambda x_2$, and $v_3 = -2\lambda x_3$
where $\lambda \neq 0$. \\
{\bf dc)} $v_1 = \lambda x_1$, $v_2 = \lambda (x_2 + x_1)$, and
$v_3 = -2 \lambda x_3$ where $\lambda \neq 0$. \\
da) There exists a cubic form $f$ such that
$y_{1,2} = \frac{\partial f}{\partial x_3}$,
$y_{2,3} = \frac{\partial f}{\partial x_1} + \lambda _2x_2x_3$,
$y_{3,1} = \frac{\partial f}{\partial x_2} - \lambda _1x_1x_3$.
The Jacobi identity then gives
$\lambda _1x_1\frac{\partial f}{\partial x_1} +
\lambda _2x_2\frac{\partial f}{\partial x_2}
+ \lambda _3x_3\frac{\partial f}{\partial x_3} = 0$,
and since $\lambda _i \neq 0$ and are pairwise different
the only possible $f = 2cx_1x_2x_3$.\\
db) There exists a cubic form $f$ such that
$y_{1,2} = \frac{\partial f}{\partial x_3}$,
$y_{2,3} = \frac{\partial f}{\partial x_1} + \lambda x_2x_3$,
$y_{3,1} = \frac{\partial f}{\partial x_2} - \lambda x_1x_3$.
The Jacobi identity then gives
$x_1\frac{\partial f}{\partial x_1}
+ x_2\frac{\partial f}{\partial x_2}
- 2x_3\frac{\partial f}{\partial x_3} = 0$,
and here
$f = g(x_1,x_2)x_3$.\\
dc) There exists a cubic form $f$ such that
$y_{1,2} = \frac{\partial f}{\partial x_3} - \frac{\lambda}{2} x_1^2$,
$y_{2,3} = \frac{\partial f}{\partial x_1} + \lambda x_2x_3$,
$y_{3,1} = \frac{\partial f}{\partial x_2} - \lambda x_1 x_3$.
The Jacobi identity then gives
$x_1 \frac{\partial f}{\partial x_1}
+ (x_1 + x_2)\frac{\partial f}{\partial x_2}
- 2x_3\frac{\partial f}{\partial x_3} = 0$,
and here $f =  cx_1^2x_3$.\\
(The easiest way to find $f$ in db) and dc) is to use the Euler identity
$x_1\frac{\partial f}{\partial x_1} +
x_2\frac{\partial f}{\partial x_2} +
x_3\frac{\partial f}{\partial x_3} = 3f$,
and to rewrite Jacobi identities accordingly.)

So in all cases a quadratic Poisson bracket up to a linear transformation
is described by a matrix in Jordan form (in case dc) it is slightly changed
for computational convenience) and by a cubic form which is
arbitrary in case a),
in variables $x_1$ and $x_2$ in case b),
$f = c_1x_1x_2x_3 + c_2x_1^3$ in case ca),
$f =  -2c_1x_3x_1^2 + c_1x_1x_2^2 + c_2x_1^3$ in case cb),
$f = 2cx_1x_2x_3$ in case da), $f = g(x_1, x_2)x_3$ in case db), and
$f =  cx_1^2x_3$ in case dc).

Now in case a) $f$ may be reduced by an arbitrary linear transformation
and in case db) by a transformation in variables $x_1, x_2$.

{\bf Remark.} It was brought to our attention that essentially the same
classification of quadratic Poisson brackets was obtained in [DH] and
[LX].

\section{Quantization of quadratic Poisson brackets in
three-dimensional case}
\label{S4}

Let $R = \k[[\hbar]]\langle a_1,a_2,a_3\rangle$ be a free algebra with
three generators
over the ring $\k[[\hbar]]$ of formal power series.
Let us define on this algebra ``Jordan'' operation $f \circ g =
\frac{1}{2}(fg + gf)$. Let $b$ be a quadratic bracket on $A = \k[x_1,x_2,x_3]$.
Then $b$ is defined by $y_{i,j} = b(x_i,x_j)$. Let $y_{i,j} =
\sum c_{i,j}^{k,l}x_kx_l$.

\begin{theorem}
\label{Th4.1}
The factor algebra $R_b$ of algebra $R$
by the ideal generated by
$[a_i,a_j] = \hbar \sum c_{i,j}^{k,l}a_k \circ a_l$
is a realization of bracket $b$. (Here $x_i$ is the image of $a_i$.)
\end{theorem}
\Proof  Let us use the same notations as above, namely
$y_{i,j} =b(x_i, x_j)$ and
$v_i = \sum_k \frac{\partial y_{i,k}}{\partial x_k}$.
Here a partial derivative is the following operation: $\circ$ is
replaced by commutative multiplication and then ordinary partial
derivative is computed.
It is easy to check that then in algebra $R$ we have
$$[x_1,y_{2,3}] + [x_2,y_{3,1}] + [x_3,y_{2,1}] =
[x_1,x_2] \circ v_3 + [x_2,x_3] \circ v_1 + [x_3,x_1] \circ v_2 \tag{1}$$
So in the case when rank is zero and all $v_i = 0$
the Theorem can be deduced from
a result of Drinfeld (see Proposition \ref{Pr2.1a}) which states that in
our setting it is sufficient to check that $W = (I_2 \otimes V) \cap
(V \otimes I_2)$ is at least one dimensional. If we denote our relations
by $f_{i,j} = y_{i,j} - [x_i, x_j]$, then\\
$f_{1,2}x_3 + f_{2,3}x_1 + f_{3,1}x_2 = x_3f_{1,2} + x_1f_{2,3} + x_2f_{3,1}$\\
which shows that $\dim(W) > 0$ in this case.
Similar proof works for the case $r = 1$ where\\
$f_{1,2}x_3 + f_{2,3}x_1 + f_{3,1}x_2 - \hbar f_{1,2}x_1 =
x_3f_{1,2} + x_1f_{2,3} + x_2f_{3,1} + \hbar x_1f_{1,2}$\\
and in the second subcase for $r = 2$ where\\
$f_{1,2}x_3 + f_{2,3}x_1 + f_{3,1}x_2 -
\hbar (c_1 f_{1,2}x_1 + f_{1,2}x_2 + f_{3,1}x_1)  =
x_3f_{1,2} + x_1f_{2,3} + x_2f_{3,1} -
\hbar (c_1 x_1f_{1,2} - x_2f_{1,2} - x_1f_{3,1})$.

Although it is clear a posteriori that $\dim(W) = 1$ in all the cases, we
approach the remaining cases differently, to avoid lengthy computations.\\

\underline{Rank two case.} In the first subcase\\
$[x_1,x_2] = \hbar c_1x_1 \circ x_2$,\\
$[x_2,x_3] = \hbar ((c_1 - \frac{\lambda}{2})x_2 \circ x_3 + 3c_2x_1^2)$,\\
$[x_3,x_1] = \hbar c_1x_3 \circ x_1$.\\

\underline{Rank three case.}

In the first subcase\\
$[x_1,x_2] = \hbar cx_1 \circ x_2$,\\
$[x_2,x_3] = \hbar (c + \frac{\lambda_2}{2})x_2 \circ x_3$,\\
$[x_3,x_1] = \hbar (c - \frac{\lambda_1}{2}) x_3 \circ x_1$.

In the second subcase it is possible to make a linear transformation of
variables $x_1$ and $x_2$ which reduces $f$ to
$(ax_1^2 + bx_1 x_2)x_3$.
So\\
$[x_1,x_2] = \frac {\hbar}{2} (2ax_1^2 + bx_1 \circ x_2)$,\\
$[x_2,x_3] = \frac {\hbar}{2} (2ax_1 \circ x_3 + (b + \lambda)x_2 \circ
x_3)$,\\
$[x_3,x_1] =  \frac {\hbar}{2} (b - \lambda) x_3 \circ x_1$.

In the third subcase relations are\\
$[x_1, x_2] = \frac {\hbar}{2} (2c - \lambda) x_1^2$,\\
$[x_2, x_3] = \frac {\hbar}{2}(2cx_1 \circ x_3 + \lambda x_2 \circ x_3)$,\\
$[x_3, x_1] = -\frac {\hbar}{2} \lambda x_3 \circ x_1$.

In all of these cases relations may be rewritten as\\
$x_2x_1 = a_3x_1x_2 + b_3x_1^2$,\\
$x_3x_2 = a_1x_2x_3 + b_1x_1^2 + c_1x_1x_3$,\\
$x_3x_1 = a_2x_1x_3$.\\

So Bergman's (see [B]) diamond lemma may be applied. It is easy to
check that in each of the cases values of parameters are such that
the only ambiguity $x_3x_2x_1$ is resolvable and the monomials
$x_1^i x_2^j x_3^k$ form a basis in the corresponding algebras.
(As it is shown in [B] resolution of this ambiguity is equivalent to
the Jacobi identity so in fact these computations show that
$\dim(W) > 0$.)

This finishes the proof of the theorem.

\section{APPENDIX: direct constructions of algebras $R_b$}
\label{S5}

In this section we give an alternative proof of the Theorem
based on direct constructions of algebras
isomorphic to $R_b$.

Let $\cal F$ be an algebra of functions on a lattice $Z^2$ with values
in $\k(\hbar, w, z)$ where $w$ and $z$ are central variables.
One of the sources of algebras isomorphic to $R_b$ will be subalgebras
of the algebra $\cal H$ of homomorphisms of $\cal F$ generated by
two ``coordinate'' shifts and by multiplication by its elements.
We will be using functional notation $f(x,y)$ for the
elements of $\cal F$ and corresponding elements of $\cal H$ and
denote shifts by  $s_x$ and $s_y$. When we need only one shift
for our construction we will denote it by $s$.
So we have commuting symbols $f(x,y)$ and relations are
$s_xf(x,y) = f(x+1,y)s_x$, $s_yf(x,y) = f(x,y+1)s_y$.

Another source will be different subalgebras and extensions of Weyl
algebras.\\

\underline{Rank zero case.}\\

In this case we can use a classification of cubic forms.

Let $H(x,y,z)$ be  a (homogeneous) cubic form over an algebraically closed
field of characteristic zero. Then the following is the list of possible
reductions of this form under the action of $GL(3)$.\\

\noindent  1) 0, 2) $x^3$, 3) $x^2y$, 4) $xy(x+y)$,\\
5) $zx^2 + xy^2$, 6) $zx^2 + y^3$,
7) $2zxy$, 8) $2zxy + x^3$, 9) $2zxy + x^3 + y^3$,\\
10) $z^2x + y(x + y)(x + cy)$ where $c \neq 0$ and $c \neq 1$.\\
A group of order 3 acts on 10).\\

This list was known of course to J. Steiner though he did not bother to
put it in quite such a form.

{\bf Remark.} First four orbits correspond to a smaller number of variables.
Forms from 5) to 9) are linear in $z$.
Form 10) is the only one with a unique singular point (at the origin).
(When $c = 0$ then 10) is equivalent to 8) and when $c = 1$ then 10) is
equivalent to 9).)\\

So let us proceed along the list of reduced  forms.\\

1) corresponds to the ring of polynomials.\\

2), 3), and 4) correspond to relations\\
$[x_1,x_2] = 0$,\\
$[x_2,x_3] = \hbar y_{2,3}$,\\
$[x_3,x_1] = \hbar y_{3,1}$,\\
where $y_{2,3}$ and $y_{3,1}$ are some elements of the
polynomial ring $P = \k[[\hbar]][x_1, x_2]$ which depend on the case.

Here we obtain skew polynomial extensions $P[x_3,\delta]$ of the
derivation type of the ring $P$ (see \cite{C}) and in fact $y_{2,3}$
and $y_{3,1}$ can be any elements of $P$. (We shall use it later.)

These algebras may be realized as subalgebras of the second Weyl
algebra $A_2$ over $\k(z)$. Say, if generators are $p_1, q_1$ and
$p_2, q_2$ (and relations are $[p_i, q_i] = 1$ and all other
commutators are zeros), then $x_1 = p_1$, $x_2 = p_2$, and
$x_3 = \hbar (q_2 y_{2,3} - q_1 y_{3,1}) + z$ will suffice.\\
Since the standard monomials in $p_1, \, p_2, \, z$ are linearly independent,
it is clear that the standard monomials in $x_1, \, x_2, \, x_3$ are linearly
independent. It implies that the the constructed algebras are isomorphic to
the corresponding $R_b$.\\

5) corresponds to relations\\
$[x_1,x_2] = \hbar x_1^2$,\\
$[x_2,x_3] = \hbar (x_1 \circ x_3 + x_2^2)$,\\
$[x_3,x_1] = \hbar x_1 \circ x_2$,\\
and  6) corresponds to relations\\
$[x_1,x_2] = \hbar x_1^2$, \\
$[x_2,x_3] = \hbar x_1\circ x_3$,\\
$[x_3,x_1] = 3\hbar x_2^2$.\\

These algebras are skew polynomial extensions $B[x_3; \alpha, \delta]$
of the algebra $B$ generated by $x_1$ and $x_2$.
(See \cite{C}, $\alpha$ is an automorphism of $B$ and $\delta$ is an
$\alpha$ derivation of $B$.)
Here $x_1^{\alpha} = x_1$,
$x_2^{\alpha} = x_2 - 2 \h x_1$ in both cases and
$x_1^{\delta} = \h x_1 \circ x_2$,
$x_2^{\delta} = - \h(x_2^2 + \h x_1 \circ x_2)$ for  5)
and $x_1^{\delta} = 3 \h x_2^2$, $x_2^{\delta} = -3 \h^2x_2^2$
for  6).

Let $A_1$ be the first Weyl algebra with generators $p$ and $q$ over
$\k(\hbar, w, z)$. (Here $w$ and $z$ are central variables and $[p,q] = 1$.)
Let $D_1$ be its field of fractions. The following subalgebras of $D_1$
are realizations of corresponding brackets.

For 5) we may consider an algebra generated by $x_1 = p^{-1}$,
$x_2 = -(\h q + w)$, and
$x_3 = -(\h q + w)^2p - \h (\h q + w) + \frac{1}{3}\h^2p^{-1} + zp^2$,
and for 6) we may consider an algebra generated by $x_1 = p^{-1}$,
$x_2 = -(\h q + w)$, and $x_3 = (\h q + w)^3p^2 + 3\h (\h q + w)^2p + zp^2$.
These algebras are homomorphic images of the corresponding $R_b$.

As above, in both these cases it is clear that the standard monomials
of constructed algebras are linearly independent (e.g., it follows from
the independence of monomials of $p^{-1}$, $w$, and $zp^2$).
So again these algebras are isomorphic to the corresponding $R_b$.\\

7), 8), and 9) correspond to relations\\
$[x_1,x_2] = \hbar x_1 \circ x_2$,\\
$[x_2,x_3] = \hbar (x_2\circ x_3 + 3c_ix_1^2)$, \\
$[x_3,x_1] = \hbar (x_3 \circ x_1 + 3d_ix_2^2)$,
where $c_7 = d_7 = d_8 = 0$ and $c_8 = c_9 = d_9 = 1$.\\

These algebras are skew polynomial extensions $Q[x_3; \alpha, \delta]$
of a quantum plane $Q$ generated by $x_1$ and $x_2$. (See \cite{C}.)
Here $x_2x_1 = kx_1x_2$ where $k = \frac{1 - \hbar}{1 + \hbar}$,
and these skew extensions are defined by
$x_1^{\alpha} = k^{-1}x_1$,
$x_2^{\alpha} = k x_2$, and
$x_1^{\delta} = \frac{c_1\hbar}{1 - \hbar}x_2^2$,
$x_2^{\delta} = \frac{c_2\hbar}{1 + \hbar}x_1^2$,
where $c_1$ and $c_2$ should be appropriately chosen for each
of the cases. It is a straightforward computation to check that
$\delta$ is well defined on $Q$ by\\
$(x_1^ix_2^j)^{\delta} =
k^{-i}(1 - k^3)^{-1}[c_3k(1 - k^{3i})x_1^{i-1}x_2^{j+2} +
c_4(1 - k^{3j})x_1^{i+2}x_2^{j-1}$\\
where $c_3 = \frac{c_1\hbar}{1 - \hbar}$,
and $c_4 = \frac{c_2\hbar}{1 + \hbar}$.

We can rewrite our relations as\\
$x_2x_1 = k x_1x_2$,\\
$x_3x_2 = k x_2x_3 + cx_1^2$,\\
$x_3x_1 = k^{-1}x_1x_3 + dx_2^2$,\\
where $k, \, c, \, d \in \k(\h)$, and
let $x_1 = s$, $x_2 = wk^{-x}s$, and $x_3 =
k(cw^{-1}k^x - dw^2k^{1-2x})(1 - k^3)^{-1}s + zk^xs^{-2}$.

As above, it is easy to check that an algebra generated by
$x_1, \, x_2$, and $x_3$ is isomorphic to the corresponding $R_b$.\\

All of the previously considered cases are also covered by Bergman's
diamond lemma (see \cite{B}).
Case 10) (which is not covered by the diamond lemma)
requires a more involved construction (see the end of the section).\\

\underline{Rank one case.}\\

Here\\
$[x_1,x_2] = 0$,\\
$[x_2,x_3] = \hbar(\frac{\partial f}{\partial x_1} - x_1x_2)$,\\
$[x_3,x_1] = \hbar\frac{\partial f}{\partial x_2}$,\\
and this case is analogous to the cases 2), 3), and 4) in
the $r = 0$ case. The corresponding algebras are skew polynomial
extensions of the ring $P$ of the derivation type and may be realized
as subalgebras of $A_2$.\\
\\

\underline{Rank two case.}\\

The first subcase where\\
$[x_1,x_2] = \hbar c_1x_1 \circ x_2$,\\
$[x_2,x_3] = \hbar ((c_1 - \frac{\lambda}{2})x_2 \circ x_3 + 3c_2x_1^2)$,\\
$[x_3,x_1] = \hbar c_1x_3 \circ x_1$\\
is analogous to the case 8) and gives a skew extension of
a quantum plane.\\
Let us rewrite relations:\\
$x_2x_1 = kx_1x_2$,\\
$x_3x_2 = k_1x_2x_3 + cx_1^2$,\\
$x_3x_1 = k^{-1}x_1x_3$.\\
The main difference with the previously considered cases is that
$k_1 \neq k$.
For this case we may take $x_1 = s_x$, $x_2 = wk^{y-x}k_1^ys_x$, and
$x_3 = cw^{-1}(1 - k_1k^2)^{-1}k^{1+x-y}k_1^{-y}s_x + zk^xs_y$.
(This representation works even when $c_1 = 0$ or
$c_1 - \frac{\lambda}{2} = 0$.)
\\

The second subcase where\\
$[x_1,x_2] = \hbar (-2c_1x_1^2$),\\
$[x_2,x_3] = \hbar ((-2c_1 + \frac {1}{2}) x_1 \circ x_3 +
(c_1 - \frac{1}{2}) x_2^2 + 3c_2x_1^2)$,\\
$[x_3,x_1] = \hbar c_1x_1 \circ x_2$\\
gives an extension of the polynomial ring $P$ when $c_1 = 0$. It is analogous
to the case 5) when $c_1 \neq 0$ (and defines an extension
of algebra $B$).

Let $c_1 = 0$. Then $x_1$ is a central element and we can take the
following two elements of the first Weyl algebra over $\k(\h, x_1, w, z)$
for $x_2$ and $x_3$ ($t$ denotes $qp$):
$x_2 = -x_1(\h t + w)$ and $x_3 = x_1(\frac{1}{2} (\h t + w)^2
- 3c_2 + zp)$.

Let us assume now that $c_1 \neq 0$.
Then we may consider algebra generated by
$x_1 = p^{-1}$, $x_2 = 2c_1(\h q + w)$, and $x_3 =
2c_1^2[(\h q + w)^2p + \h (\h q + w)] +
(c_1^2\h^2 (4c_1 - 1) - 3c_2 )(1 - 6c_1)^{-1}p^{-1} + zp^d$
where $d = (2c_1)^{-1} - 2$. If $d$ is not an integer we can view this
algebra as a subalgebra of $D_1[p^d]$.
When $1 - 6c_1 = 0$ this algebra is not defined and we put
$x_3 =  2c_1^2[(\h q + w)^2p + \h (\h q + w)] +
(c_1^2\h^2 (4c_1 - 1) - 3c_2 )(2c_1)^{-1}p^{-1}\ln(p) + zp^{-1}$.
This algebra belongs to $D_1[\ln(p)]$.

These two extensions of $D_1$ are defined as follows.
Each element $a \in D_1$ may be represented as
$a = \sum^{i=k}_{-\infty} a_i p^i$ where $a_i \in \k(\h, w, z, q)$.
It is possible to extend multiplication from $D_1$ on such arbitrary series
by $a \star b = \sum_{j=0}^{\infty} \frac{1}{j!}
\frac{\partial^j a}{\partial p^j} \frac{\partial b}{\partial q^j}$,
where  multiplication in the right side of this formula is ordinary
commutative multiplication of power series.

It is clear that this formula gives well defined expressions when either
$a$ or $b$ are replaced by $f$ if $\frac{\partial f}{\partial p}$ and
$\frac{\partial f}{\partial q}$ are defined and are ''sufficiently small".
For example, in our cases when $f = p^d$ or $f = \ln(p)$ we have
$\frac{\partial f}{\partial p} = p^{-1}(af + b)$ and
$\frac{\partial f}{\partial q} = 0$ and $\star$ multiplication is well
defined.
\\

\underline{Rank three case.}\\

The first subcase where\\
$[x_1,x_2] = \hbar cx_1 \circ x_2$,\\
$[x_2,x_3] = \hbar (c + \frac{\lambda_2}{2})x_2 \circ x_3$,\\
$[x_3,x_1] = \hbar (c - \frac{\lambda_1}{2}) x_3 \circ x_1$\\
is analogous to the case 7) and is a ``quantum space''.\\
Let us rewrite these relations as
$x_2x_1 = k_3 x_1x_2$,
$x_3x_2 = k_1 x_2x_3$,
$x_3x_1 = k_2x_1x_3$
where $k_i \in \k (\h )$.
We may put $x_1 = s_x$, $x_2 = k_3^{-x}s_y$, and $x_3 = zk_2^{-x}k_1^{-y}$.
\\

In the second subcase it is possible to make a transformation in
variables $x_1$ and $x_2$ and to reduce $f$ to either\\
a) $f = 0$ or\\
b) $ f = x_1x_2x_3$ or\\
c) $f = x_1^2x_3$.\\
For a) relations are\\
$[x_1,x_2] = 0$,\\
$[x_2,x_3] = \hbar \frac{\lambda}{2}x_2 \circ x_3$,\\
$[x_3,x_1] = - \hbar  \frac{\lambda}{2} x_3 \circ x_1$,\\
and the corresponding algebra may be considered
as a quantum space.\\
For b) relations are\\
$[x_1,x_2] = \frac{\hbar}{2} x_1 \circ x_2$,\\
$[x_2,x_3] = \hbar  \frac{\lambda + 1}{2}x_2 \circ x_3$,\\
$[x_3,x_1] =  \hbar  \frac{1 - \lambda}{2} x_3 \circ x_1$,\\
and this is a quantum space again.\\
For c) relations are\\
$[x_1,x_2] = \hbar x_1^2$,\\
$[x_2,x_3] = \hbar (x_1 \circ x_3 + \frac{\lambda}{2}x_2 \circ x_3)$,\\
$[x_3,x_1] = - \hbar  \frac{\lambda}{2} x_3 \circ x_1$.\\

In the third subcase relations are\\
$[x_1, x_2] = \frac {\hbar}{2} (2c - \lambda) x_1^2$,\\
$[x_2, x_3] = \frac {\hbar}{2}(2cx_1 \circ x_3 + \lambda x_2 \circ x_3)$,\\
$[x_3, x_1] = -\frac {\hbar}{2} \lambda x_3 \circ x_1$.\\

Here is a realization for these two cases.
Let us consider the following relations:\\
$[x_1,x_2] = ax_1^2$,\\
$[x_2,x_3] = (bx_1 - cx_2) \circ x_3$,\\
$[x_3, x_1] = cx_1 \circ x_3$.\\
Let us rewrite the last relation as $x_3x_1 = kx_1x_3$.\\
Then $x_3x_2 = (kx_2 + k_1x_1)x_3$, and we have a skew
polynomial extension of automorphism type of the subalgebra
generated by $x_1, \, x_2$.

If $a = 0$, then it is an extension of polynomial algebra which can be
realized as a subalgebra of $\cal H$ with $x_1 = wk^x$,
$x_2 = wk_1xk^{x-1}$, and $x_3 = s$.

If $a \neq 0$ then a realization can be found in the tensor product of
$D_1[p^d]$ and $\cal H$ (which is rather natural since it contains both
a quantum plane and a big subalgebra of a Weyl algebra). We can take
$x_1 = p^{-1}s$, $x_2 = -(aq + w)s$ and $x_3 = zp^dk^{-x}$ where
$d = -k_1k^{-1}a^{-1}$.
\\

\underline{Orbit 10) case.}\\

To construct algebras which correspond to the case 10),
one may consider an algebra $B =  \k(u)[[\hbar, p, q]]$
of formal power series in $\h$, $p$, and $q$ over $\k(u)$
(where $u$ is a central variable and $[p,q] = \h$) and search
for a realization
of corresponding brackets as subalgebras of $B$ with generators
$x_1 = \sum \hbar^iq^jx_{1,i,j}$,
$x_2 = \sum \hbar^iq^jx_{2,i,j}$,
and $x_3 = p$.
Here $x_{1,i,j}$, $x_{2,i,j} \in \k(u)[[p]]$.

With the help of rather straightforward computations one may find
corresponding solutions and show that the standard monomials
are linearly independent in this case as well.\\

Let us assume that $H(x,y,z) = z^2x + c_1xy^2 + c_2x^2y + y^3$.
Then\\
$[x_1,x_2] = \hbar x_1 \circ x_3$,\\
$[x_2,x_3] = \hbar (x_3^2 + c_1x_2^2 + c_2 x_1 \circ x_2)$, \\
$[x_3,x_1] = \hbar (c_1x_1 \circ x_2 + c_2x_1^2 + 3x_2^2)$\\
where we perceive these relations as equations for $x_{1,i,j}$
and $x_{2,i,j}$.\\
Let us again use notations $[x_i,x_j] = \hbar y_{i,j}$ and rewrite
the second and third equations as\\
$\sum j \hbar^{i+1}q^{j-1}x_{2,i,j}  = - \hbar y_{2,3}$ and\\
$\sum j \hbar^{i+1}q^{j-1}x_{1,i,j}  = \hbar y_{3,1}$.\\
Now $g(p)q^i = \sum_{k=0}^i \hbar^k \frac{i!}{k!(i-k)!} q^{i-k}g^{(k)}$
where $g(p) \in \k[[p]]$ and $g^{(k)}$ is the ``ordinary"
derivative.
So the product of two monomials in $B$ is a linear combination of
monomials with the same total $\h, q$-degree (which is equal to the
sum of degrees of multiples) and $\h$-degrees of these monomials are
not less than the sum of $\h$ degrees of the multiples.
Therefore $x_{m,i,j}$ can be expressed through $x_{1,a,b}$ and $x_{2,a,b}$
where $a + b = i + j - 1$ and $a \leq i$.
So by induction one can show that $x_{m,i,j}$ can be expressed
as polynomials of $x_{1,n,0}^{(k)}$ and $x_{2,n,0}^{(k)}$ where $n \leq i$.
Since $x_{m,n,0}$ do not appear in the left hand sides of our equations they
can be chosen arbitrarily.\\
Now we have to choose these coefficients in such a way that the first
equation will be satisfied.
The following observation make this part of the computations
more manageable.
It is easy to check that
$[x_1,y_{2,3}] + [x_2,y_{3,1}] + [x_3,y_{2,1}] =
[x_1,x_2] \circ v_3 + [x_2,x_3] \circ v_1 + [x_3,x_1] \circ v_2 = 0$
for any choice of $x_i$.
On the other hand, with our choice of $x_i$ the Jacobi identity gives
$\h[x_1,y_{2,3}] + \h[x_2,y_{3,1}] + [x_3,[x_1, x_2]] = 0$.
It implies that $[p,[x_1,x_2] - \hbar x_1 \circ p] = 0$
is satisfied with any choice of $x_{m,n,0}$.
Therefore, $[x_1,x_2] - \hbar x_1 \circ p$ does not contain terms with $q$,
and it is sufficient to check that $[x_1,x_2] = \hbar x_1 \circ x_3$
is an equality only in terms which do not contain $q$.
Let us use notation $\equiv$ for such equalities.\\
It is clear then that
$[x_1,x_2] \equiv
[ \sum \hbar^ix_{1,i,0}, x_2] +
[x_1, \sum \hbar^i x_{2,i,0}]$,\\
and that
$p \circ x_1 \equiv
2 \sum \hbar^i x_{1,i,0}p  +
\hbar \sum \hbar^ix_{1,i,1}$.\\
So the coefficient $\Delta_n$ of $\hbar^n$ in
$[x_1,x_2] - \hbar x_1 \circ p$ when $n > 1$
is equal to \\
$x_{1,n-1,0}^{(1)}x_{2,0,1} -
x_{2,n-1,0}^{(1)}x_{1,0,1} +
x_{1,0,0}^{(1)}x_{2,n-1,1} -
x_{2,0,0}^{(1)}x_{1,n-1,1} -
2x_{1,n-1,0}p + \delta_n$\\
where $\delta_n$ depends on $x_{m,k,0}$ with $k < n-1$.\\
The coefficient $\Delta_1$ of $\hbar$ is
$x_{1,0,0}^{(1)}x_{2,0,1} -
x_{2,0,0}^{(1)}x_{1,0,1} -
2x_{1,0,0}p$.\\
Now\\
$x_{1,k,1} =
\frac{\partial^2 H}{\partial x_1 \partial x_2}(x_{1,0,0},x_{2,0,0},p)x_{1,k,0}
+ \frac{\partial^2 H}{\partial x_2 \partial
x_2}(x_{1,0,0},x_{2,0,0},p)x_{2,k,0}
+ \delta_{1,k}$\\
and\\
$-x_{2,k,1} =
\frac{\partial^2 H}{\partial x_1 \partial x_1}(x_{1,0,0},x_{2,0,0},p)x_{1,k,0}
+ \frac{\partial^2 H}{\partial x_1 \partial
x_2}(x_{1,0,0},x_{2,0,0},p)x_{2,k,0}
+ \delta_{2,k}$\\
for $k > 0$.
(Here $\delta_{m,k}$ depend on $x_{m,i,0}$ with $i < k$.)\\
For $k = 0$ one gets\\
$x_{1,0,1} =
\frac{\partial H}{\partial x_2}(x_{1,0,0},x_{2,0,0},p)$\\
and\\
$-x_{2,0,1} =
 \frac{\partial H}{\partial x_1}(x_{1,0,0},x_{2,0,0},p)$.\\
By substituting these expressions in $\Delta_1$
one gets\\
$ \frac{\partial H}{\partial
x_1}(x_{1,0,0},x_{2,0,0},p)x^{(1)}_{1,0,0} +
\frac{\partial H}{\partial
x_2}(x_{1,0,0},x_{2,0,0},p)x^{(1)}_{2,0,0} + 2x_{1,0,0}p = 0$.\\
Since $2x_{1,0,0}p = \frac{\partial H}{\partial
x_3}(x_{1,0,0},x_{2,0,0},p)$ one can conclude that
$\frac{\partial H}{\partial p} = 0$ and put
$H(x_{1,0,0},x_{2,0,0},p) = u$ where $u$ is a central variable.\\
Substitution in $\Delta_n$ for $n > 1$ gives\\
$  \frac{\partial}{\partial p}(\frac{\partial H}{\partial
x_1}(x_{1,0,0},x_{2,0,0},p)x_{1,n-1,0} +
\frac{\partial H}{\partial
x_2}(x_{1,0,0},x_{2,0,0},p)x_{2,n-1,0}) + \epsilon_n = 0$\\
accordingly.
(Here $\epsilon_n$ depend on $x_{m,i,0}$ with $i < n-1$.)\\
Let us put finally $x_{1,i,0} = 0$ for all $i$ and $x_{2,0,0}
= v$ where $v$ is a central variable. Then $u = v^3$ and
$x_1 = 3v^2q + ...$, $x_2 = v + ...$, $ x_3 = p$ where omitted terms are of
higher  $\hbar, q$-degree. Since the standard monomials in the lowest terms
are linearly independent, we can conclude that the standard monomials in
$x_1, x_2, x_3$ are also linearly independent.

{\bf Remark.} Algebras which appeared here are considered
from a different point of view in \cite{AS}, \cite{ATV}, \cite{FO1},
\cite{FO2}, \cite{O}.


\small

\begin{thebibliography}{ATV}
\bibitem[AS]{AS} M. Artin and W. Schelter,
{\em Graded Algebras of Global Dimension 3\/},
Advances in Math. 66 (1987), 171-216.
\bibitem[ATV]{ATV} M. Artin, J. Tate, M. Van Den Bergh,
{\em Some Algebras Associated to Automorphisms of Elliptic Curves\/},
The Grothendieck Festschrift, v. 1, Birkhauser, 1990.
\bibitem[Ba]{Ba} J. Bakelin,
{\em A distributiveness property of augmented algebras and some
related homological results\/},
Preprint, Univ. Stockholm, 1981.
\bibitem[B]{B} G.M. Bergman, {\em The Diamond Lemma for Ring
Theory\/},
 Advances in Mathematics, 29 (1978), 178-218.
\bibitem[C]{C} P.M. Cohn, ``Skew Field Constructions'',
 Cambridge University Press, 1977.
\bibitem[D]{D} V. Drinfeld,
{\em On quadratic commutator relations in the quasiclassical case\/},
Math. physics and functional analysis, Kiev, Naukova Dumka, 1986.
\bibitem[DH]{DH} J. Dufour et A. Haraki,
{\em Rotationnels et structures de Poisson quadratiques\/},
C. R. Acad. Sci. Paris, V. 312 I (1991), 137-140.
\bibitem[FO1]{FO1} B. Feigin and A. Odessky,
{\em Sklyanin Algebras Associated with Elliptic Curves\/},
Preprint, Institute for Theoretical Physics, Ukraine, Kiev, 1989.
\bibitem[FO2]{FO2} B. Feigin and A. Odessky,
{\em Elliptic Sklyanin Algebras\/}, Functional Analysis and Applications,
V. 23, 3, 1989.
\bibitem[LX]{LX} Z. Liu and P. Xu,
{\em On quadratic Poisson structures\/},
Letters in Math. Phys., 26 (1992), 33-42.
\bibitem[O]{O} A. Odessky,
{\em Infinite  Dimensional Algebras and Complex Varieties\/},
Doctoral Thesis, Moscow State University, Moscow, 1991.
\bibitem[PP]{PP} A. Polistchuk and L. Posicelsky,
{\em On Quadratic Algebras\/}, Preprint, Moscow, 1991.
\end{thebibliography}
\end{document}